\documentclass[12pt,reqno]{amsart}

\usepackage{amsfonts,amsmath,amsthm}
\usepackage{amssymb,epsfig,color}
\usepackage{enumerate} 
\usepackage{mathrsfs}

\usepackage{color} 
\definecolor{vert}{rgb}{0,0.6,0}

\usepackage[text={425pt,650pt},centering]{geometry}
\usepackage{graphicx}
\usepackage{epsfig}
\usepackage{tikz}
\usepackage{caption}
\usepackage{color} 
\definecolor{vert}{rgb}{0,0.6,0}

\numberwithin{figure}{section}

\theoremstyle{plain}
\newtheorem{thm}{Theorem}[section]

\newtheorem{quest}{Question}

\newtheorem{lem}[thm]{Lemma}
\newtheorem{cor}[thm]{Corollary}

\theoremstyle{remark}
\newtheorem{rem}{\bf{Remark}}
\numberwithin{equation}{section}

\newcommand{\N}{\mathbb{N}}

\newcommand{\R}{\mathbb{R}}

\newcommand{\T}{\mathbb{T}}
\newcommand{\Z}{\mathbb{Z}}

\newcommand{\BUC}{{\rm BUC\,}}

\newcommand{\Lip}{{\rm Lip\,}}

\newcommand{\al}{\alpha}
\newcommand{\gam}{\gamma}
\newcommand{\del}{\delta}
\newcommand{\ep}{\varepsilon}

\newcommand{\Gam}{\Gamma}

\newcommand{\ol}{\overline}

\begin{document}

\title[Effective fronts of polygon shapes in two dimensions]
{Effective fronts of polygon shapes in two dimensions}

\author[W. JING, H. V. TRAN, Y. YU]
{Wenjia Jing, Hung V. Tran, Yifeng Yu}

\thanks{
The work of WJ is partially supported by the NSFC Grant No.\,11871300.
The work of HT is partially supported by NSF CAREER grant DMS-1843320 and a Simons Fellowship.
The work of YY is partially supported by NSF grant DMS-2000191.
}

\address[W. Jing]
{
Yau Mathematical Sciences Center,
Tsinghua University, No.1 Tsinghua Yuan,
Beijing 100084, China }
\email{wjjing@tsinghua.edu.cn}

\address[H. V. Tran]
{
Department of Mathematics, 
University of Wisconsin Madison, 480 Lincoln  Drive, Madison, WI 53706, USA}
\email{hung@math.wisc.edu}

\address[Y. Yu]
{
Department of Mathematics, 
University of California, Irvine,
410G Rowland Hall, Irvine, CA 92697, USA}
\email{yyu1@math.uci.edu}

\date{\today}

\maketitle

\begin{abstract}
We study the effective fronts of first order front propagations in two dimensions ($n=2$) in the periodic setting.
Using PDE-based approaches, we show that for every $\alpha\in (0,1)$,  the class of centrally symmetric polygons with rational vertices and nonempty interior is admissible as effective fronts for given front speeds in $C^{1,\alpha}(\T^2,(0,\infty))$. 
This result can also be formulated in the language of stable norms corresponding to periodic metrics in $\T^2$.
Similar results were known long time ago when $n\geq 3$  for front speeds in $C^{\infty}(\T^n,(0,\infty))$.  
Due to topological restrictions, the  two dimensional case is much more subtle. 
In  fact,  the effective front is $C^1$, which cannot be a polygon, for given  $C^{1,1}(\T^2,(0,\infty))$ front speeds.
Our regularity requirements on front speeds  are hence optimal. 
To the best of our  knowledge, this is the first time that polygonal effective fronts have been constructed in two dimensions. 
\end{abstract}


\section{Introduction}

In this paper, we are concerned with fine properties of the effective fronts of first order front propagations in oscillatory periodic environment in two dimensions.
The front propagation is modeled by the following Hamilton-Jacobi equation with oscillatory periodic coefficient:
\begin{equation}\label{HJ-ep}
\begin{cases}
u^\ep_t + a\left(\frac{x}{\ep}\right) |Du^\ep| = 0 \quad &\text{ in } \R^n \times (0,\infty),\\
u^\ep(x,0) = g(x) \quad &\text{ on } \R^n.
\end{cases}
\end{equation}
Here, $\ep \in (0,1)$, $g\in \BUC(\R^n) \cap \Lip(\R^n)$ is the initial data, where $\BUC(\R^n)$ is the space of bounded, uniformly continuous functions on $\R^n$. The coefficient $a: \R^n \to (0,\infty)$ determines the normal velocity in the underlying front propagation model. 
Throughout the paper we deal with $a$ that is continuous,  $\Z^n$-periodic and non-constant positive function.
Denote by $\T^n=\R^n/\Z^n$ the $n$-dimensional flat torus.
We then write $a \in C(\T^n,(0,\infty))$.

\smallskip

We now give a minimalistic review of the literature on the qualitative homogenization of \eqref{HJ-ep}, which fits in the classical and standard framework (see \cite{LPV, Ev1,Tr}). 
As $\ep \to 0$, the solution $u^\ep$ of \eqref{HJ-ep} converges, locally uniformly on $\R^n \times [0, \infty)$, to the solution of the effective (homogenized) problem:
\begin{equation}\label{HJ}
\begin{cases}
u_t + \ol{H}(Du) = 0 \quad &\text{ in } \R^n \times (0,\infty),\\
u(x,0) = g(x) \quad &\text{ on } \R^n.
\end{cases}
\end{equation}
Here, $\ol{H}$ is the effective Hamiltonian, which is determined by the Hamiltonian $H(x,p) =a(x)|p|$ of \eqref{HJ-ep} in a nonlinear way as follows. 
For each $p\in \R^n$, $\ol{H}(p)$ is the unique real  number such that the following equation admits a continuous viscosity solution
\[
 {\rm (E)}_p \qquad a(y) |p +Dv_p(y)| = \ol{H}(p) \quad \text{ in } \T^n.
\] 
This is the well-known cell (ergodic) problem. 
Although $\ol H(p)$ is unique, $v_p$ is not unique in general even up to additive constants.
 It is known that $\ol H(p)$ has the following inf-max representation formula (see e.g., \cite{Tr})
\begin{align}\label{inf-max}
\ol H(p)=\inf_{\phi\in C^{\infty}(\Bbb T^n)}\max_{y \in \T^n}\,a(y)|p+D\phi(y)|
=\inf_{\phi\in C^{1}(\Bbb T^n)}\max_{y \in \T^n}\,a(y)|p+D\phi(y)|.
\end{align}
Clearly,  $\overline H$ is convex, even, and positively homogeneous of degree $1$. 
We sometime write $\overline H=\overline H_a$ to emphasize the dependence on the function $a$.
Due to those properties of $\overline H_a$, its $1$-sublevel set
\[
S_{a} :=\{p\in \R^n \,:\,  \overline H_a(p) \leq 1\} 
\]
belongs to $\mathcal{W}$, which denotes the collection of all convex sets in $\R^n$ that are centrally symmetric with nonempty interior. 

Due to the 1-homogeneity of $\ol H_a$, it is determined by $S_a$. By \emph{effective front} we mean the convex set $D_a$ dual to $S_a$ in $\R^n$. 
It is known that $D_a$ is the large time averaged limit of the so-called reachable set that arises in the control representation of \eqref{HJ-ep}.
See more discussions in Remark \ref{stablenorm}.

\smallskip

As in general homogenization theory, the function $a$ in \eqref{HJ-ep} models the periodic environment that rules the front propagation in the microscopic scale. In the limit of $\ep \to 0$, the homogenized problem \eqref{HJ} captures the effects of the oscillatory periodic environment on front propagations in the macroscopic level.
From both mathematical and practical point of views, it is very important and interesting to characterize 
certain finer details of the effective Hamiltonian $\overline H$, or equivalently, those for the set $S_a$ in the current setting. 
For example,  in combustion literature,  the well-known G-equation is often used as another front propagation model, and the effective burning velocity associated to it is sometimes taken to be isotropic for convenience (see \cite{KTW}).   
This strongly motivates the following question.

\begin{quest} \label{Q1}
For what kind of $W\in  \mathcal{W}$ does there exist a function $a \in C(\T^n,(0,\infty))$ such that $S_a=W$? 
\end{quest}

Here, we set $a$ in the regularity class $C(\T^n,(0,\infty))$ since this is most common in the homogenization theory of Hamilton-Jacobi equations. 
In  environments of composite materials,  piecewise constant functions (or more generally, $L^{\infty}(\T^n,(0,\infty))$ functions) are probably more suitable.

The above question  is often called the realization problem, which remains largely open. 
Below we summarize what is known so far. 
Most of them were formulated in the equivalent forms in terms of the $\beta$ functions in Aubry-Mather theory or in terms of the stable norms of periodic metrics in geometry; see Remark \ref{stablenorm}. 

\begin{enumerate}
\item When $n\geq 3$,  all centrally symmetric polytopes with rational vertices and  nonempty interior are realizable for some front speed $a\in C^{\infty}(\T^n, (0,\infty))$.   
This was first studied in the seminal work \cite{Hed}, and completed by \cite{Ban,BB2006,J2009,JTY2} via different approaches. 
Hence,  the set of realizable convex sets are dense in  $\mathcal{W}$.  
Very little is known about finer properties of $S_a$  except along some irrational directions (see \cite{BIK} for instance).  

\item When $n=2$,   the result in   \cite{Car} implies  that $S_a$ is $C^1$  if $a\in C^{1,1}(\Bbb T^2)$ due to the two dimensional topological restriction and the fact that the initial value problem of the ODE system $\dot \xi=V(\xi)$ has a unique solution for Lipschitz continuous $V$. 
 If we assume further that $a\in C^{\infty}(\T^2, (0,\infty))$, then there are other restrictions on $S_a$.  For example, \cite{Ban2} yields that $S_a$ cannot be a strictly convex set other than a disk. 
\end{enumerate}

A very natural  question arising from (2) above is: for $n=2$, are all centrally symmetric polygons with rational slopes and nonempty interior realizable if we lower the regularity of $a$ to $C^{1,\alpha}(\T^2, (0,\infty))$ for $\alpha \in (0,1)$?
See \cite{Burago,Bol} for more questions on possible shapes and differentiability properties of $S_a$.
Hereafter,  a polygon is said to  be {\it centrally symmetric  with rational slopes}  if it can be expressed as 
\begin{equation}
\label{eq:csprs}
P=\{p\in \R^n \,:\,   \max_{1 \leq i \leq m}|q_i\cdot p| \leq 1\}
\end{equation}
for $m$ given rational vectors $\{q_i\}_{i= 1}^{m} \subset \R^n$. 

\medskip

The following is our main result, which gives an affirmative answer to the above question.

\begin{thm} \label{thm:main}
Assume that $n=2$.
Then,  for any $\alpha\in (0,1)$ and for any centrally symmetric polygon $P$ with rational slopes and nonempty interior, there exists $a \in C^{1, \alpha}(\Bbb T^2, (0,\infty))$ such that 
\[
S_a=P.
\]
\end{thm}

\begin{rem} 
The functions $a=a(x)$ in our constructions are also $C^2$ away from finitely many points on $\Bbb T^2$. 
It is not hard to show that $S_a$ cannot have an edge of irrational slope within this class of functions.
\end{rem}

\begin{rem}\label{stablenorm}
We can formulate the result of Theorem \ref{thm:main} in the language of stable norms corresponding to periodic metrics as follows. 
We view $a\in C(\mathbb{T}^2;(0,\infty))$ as a $\mathbb{Z}^2$-periodic function in $\R^2$, and it defines a Riemannian metric 
\[
g =\frac{1}{a(x)}(dx_{1}^2+dx_{2}^2)
\]  
on $\R^2$ that is clearly periodic. 
Let $d_a(\cdot,\cdot)$ denote the distance function induced by this metric.
The \emph{stable norm} associated to $g$, or $a$, is well defined by
\[
\|v\|_{a}=\lim_{\lambda\to \infty}{d_a(0,\lambda v)\over \lambda}, \qquad v\in \R^2.
\]
See \cite{Burago} for more background. 
In particular, it was proved there that, for all $v\in \R^2$,
\[
\left|\lambda \|v\|_a-d_a(0,\lambda v)\right|\leq C
\]
for a universal  constant $C>0$. 
For $P$ with rational slopes $\{q_i\}_{i= 1}^{m}$, $S_a=P$ means that the unit ball $\ol B_1^a$ of the stable norm $\|\cdot\|_a$ satisfies
\[
\ol B_1^a= {\rm conv}\left(\{\pm q_i\,:\, 1\leq i \leq m\}\right).
\]
Here, ${\rm conv}(E)$ is the convex hull of a set $E$. 
In fact, $\ol B_1^a$ is the dual convex set of $S_a$ in $\R^n$.
Thus, Theorem \ref{thm:main} implies that all centrally symmetric polygons with rational vertices and nonempty interior are obtainable as unit balls of the stable norms corresponding to some $a  \in C^{1, \alpha}(\Bbb T^2, (0,\infty))$. 
In \cite{JTY2}, we wrote $\ol B_1^a$ as $D_a$.

In view of the duality between $S_a$ and $\ol B^a_1$, Question \ref{Q1} can also be reformulated in terms of unit balls of stable norms.

\begin{quest} \label{Q2}
For what kind of $W\in  \mathcal{W}$ does there exist a function $a \in C(\T^n,(0,\infty))$ such that $\ol B_1^a=W$?
\end{quest}
\end{rem}

We mention that it was proved in \cite{MPS}  that there exists $a\in C^{\infty}(\T^2, (0,\infty))$ such that 
\begin{equation}\label{partial}
\{\pm q_i\,:\, 1\leq i \leq m\}\subset \partial \ol B_1^a,
\end{equation}
i.e., the  stable norm is partially prescribed. 
Our PDE based approach also provides a very simple proof of this fact. 
See Remark \ref{newproof} at the end of Section \ref{sec:main}. 

\medskip

It is worth mentioning that Questions \ref{Q1}--\ref{Q2} also appear in the first passage percolation community, where the unit ball of the stable norm is called the limit shape.
In the general stationary ergodic setting, that is $a:\R^n \times \Omega \to (0,\infty)$ is stationary ergodic, it was proved in \cite{CoxDur} that the limit shape exists and is a deterministic convex compact set in $\R^n$.
Then, it was shown in \cite{HM} that any 
symmetric compact convex set $C$ with nonempty interior 
is a limit shape corresponding to some stationary ergodic $a$.
However, when $a$ is restricted to the independent and identically distributed (i.i.d.) setting, Question \ref{Q2} is completely open, and it is of great interests to study properties of the limit shape.
We refer the readers to \cite{ADH} for detailed discussions and a list of extremely interesting open problems.
For example, it is expected that the $n$-dimensional cube is not a possible limit shape in the i.i.d. setting.

\medskip

As an immediate consequence of Theorem \ref{thm:main}, we  obtain the following result,  which also follows from  the less delicate inclusion \eqref{partial}. 
\begin{cor}
The two collections
\[
\left\{S_a\,:\, a \in C^\infty(\T^2, (0,\infty))\right\} \quad \text{ and } \quad \left\{\ol B_1^a \,:\, a \in C^\infty(\T^2, (0,\infty))\right\}
\]
are both proper dense sets in $\mathcal{W}$. 
\end{cor}

Our proof of Theorem \ref{thm:main} is done by explicit construction and relies on the characterization \eqref{inf-max}. Similar to that of the higher dimension cases,  a rough idea to construct $a$ is sort of clear:  form a network of curves pointing to the rational directions $\{q_i\}_{i= 1}^{m}$, and assign values of $a$ appropriately in this network. The curves in this network serve as highways so that proper assignment of values of $a$ here guarantees that $q_i$'s are in the effective front, proving the lower bound. Let $a$ be very small away from the network of highways so that, to check the upper bound, we can choose appropriate test functions in \eqref{inf-max} and still concentrate on behaviors close to the network.  
In three dimensions, this strategy is easy to be carried out by just choosing disjoint straight lines pointing to the directions $\{q_i\}_{i= 1}^{m}$, thanks to the availability of space.  
In two dimensions, however, those curves always intersect, and it is very delicate  to design $a$ near the intersection points of the network and to make everything compatible.  

We would like to mention that this paper belongs to the ongoing project of systematic studies of inverse problems in periodic homogenization of Hamilton-Jacobi equations 
(see \cite{LTY,JTY,TY}).

\subsection*{Outline of the paper}
The proof of Theorem \ref{thm:main} is given in Section \ref{sec:main}.
Some auxiliary results are given in Appendix \ref{appen}.

\section{Proof of Theorem \ref{thm:main}}\label{sec:main}


Let $P$ be a centrally symmetric polygon with rational slopes $\{q_i\}_{i=1}^m$ of the form \eqref{eq:csprs}. Since $n=2$, we can assume that the rational vectors $\{q_i\}_{i=1}^{m} \subset \R^2$ are arranged clockwise; see Figure \ref{fig1}.
For each $i=1,..,m$, there is a unique real number $\lambda_i > 0$, and a unique irreducible integer vector $(m_i,n_i) \in \Z^2$ so that
\[
q_i=\lambda_i\,(m_i,n_i).
\]
Note that by the definition \eqref{eq:csprs}, $\{q_i\}_{i=1}^m$ form normal vectors of half of the faces of $P$. By symmetry, we order the other half by 
\[q_{m+i}=-q_i, \qquad 1\leq i \leq m.
\]
Let $p_i$ be the vertex between $q_i$ and $q_{i+1}$ for $1 \leq i \leq 2m-1$. 
Then the vertices $\{p_i\}_{i=1}^m$ of $P$ (in fact, half of them) are determined by
\begin{equation}\label{strict}
p_i\cdot q_i=p_i\cdot q_{i+1}=1 \quad \mathrm{ and } \quad \max_{\substack{j\not= i,i+1\\ 1\leq j \leq m}}|q_j\cdot p_i|<1, \qquad 1\le i\le m.
\end{equation}

\begin{center}
\includegraphics[scale=0.6]{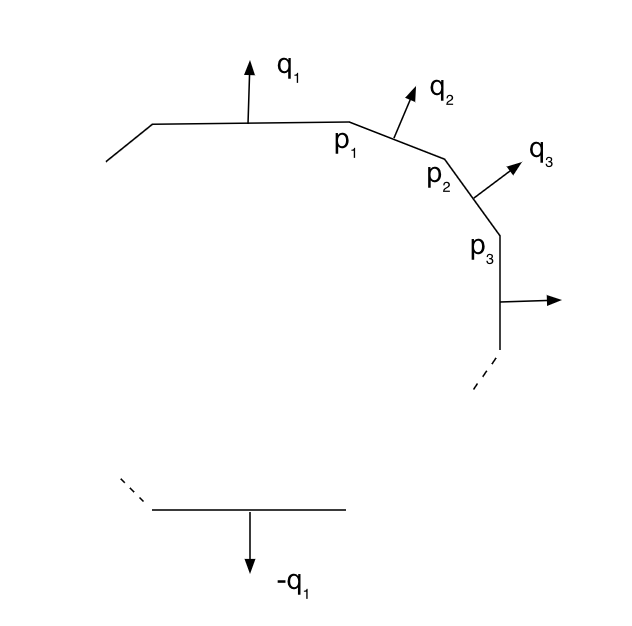}
\captionof{figure}{Polygon $P$ with vertices $p_1, p_2, \ldots, p_{2m}$}\label{fig1}
\end{center}

\begin{lem}\label{one-sided} 
Suppose that $\xi\in C^1([0,T], \R^2)$ satisfies that 
\[
\xi (T)-\xi (0)=(m,n)\in \Z^2.
\]
Then
\[
\overline H(p)\geq \lambda \,  p\cdot (m,n)
\]
for ${1\over \lambda} :=\int_{0}^{T} {1\over a(\xi(t))}|\dot \xi(t)|\,dt$. 
\end{lem} 

\begin{proof} 
Owing to the inf-max formula \eqref{inf-max}, it suffices to show that for any $\phi\in C^{\infty}(\Bbb T^2)$,
\[
M :=\max_{x\in \R^2} a(x)|p+D\phi(x)|\geq \lambda \, p\cdot (m,n).
\]
Let $u(x)=p\cdot x+\phi(x)$. 
Then, in view of the computation
\[
p\cdot (m,n)=u(\xi(T))-u(\xi(0))=\int_{0}^{T} Du(\xi(t))\cdot \dot \xi(t)\,dt\leq {M\over \lambda},
\]
the desired inequality follows immediately.
\end{proof}

\medskip

\begin{proof}[Proof of Theorem \ref{thm:main}]
We divide the proofs into several steps.

\medskip

\noindent {\bf Step 1.  Creation of a suitable network.} 
First we  choose $m$ lines $\{L_i\}_{i=1}^{m}$ in $\R^2$ such that $L_i$ is parallel to $q_i$ and, when projected to $\Bbb T^2$,  no three lines intersect at one point. 
Then, by \eqref{strict},  for every two distinct points $x$ and $y$ on $L_i$, we have that 
\[
|p_i\cdot (x-y)|> \max_{\substack{j\not =i-1, i\\ 1\leq j \leq m}}|p_j\cdot (x-y)|.
\]
Consider all  integer translations of $L_i$, which form a network  $\cup_{i=1}^{m}\left(L_i+\Z^2\right)$.  
Let
\[
\text{$I=$ the collection of all intersection points in $\bigcup_{i=1}^{m}\left(L_i+\Z^2\right)$}.
\]
Note that the intersection set $I$ is $\Z^2$-periodic. 
Denote 
\[
d=\min\{|x-y|:\   x\not=y, \  x, y\in I\}.
\]
Clearly, $d>0$. 

\begin{center}
\includegraphics[scale=0.5]{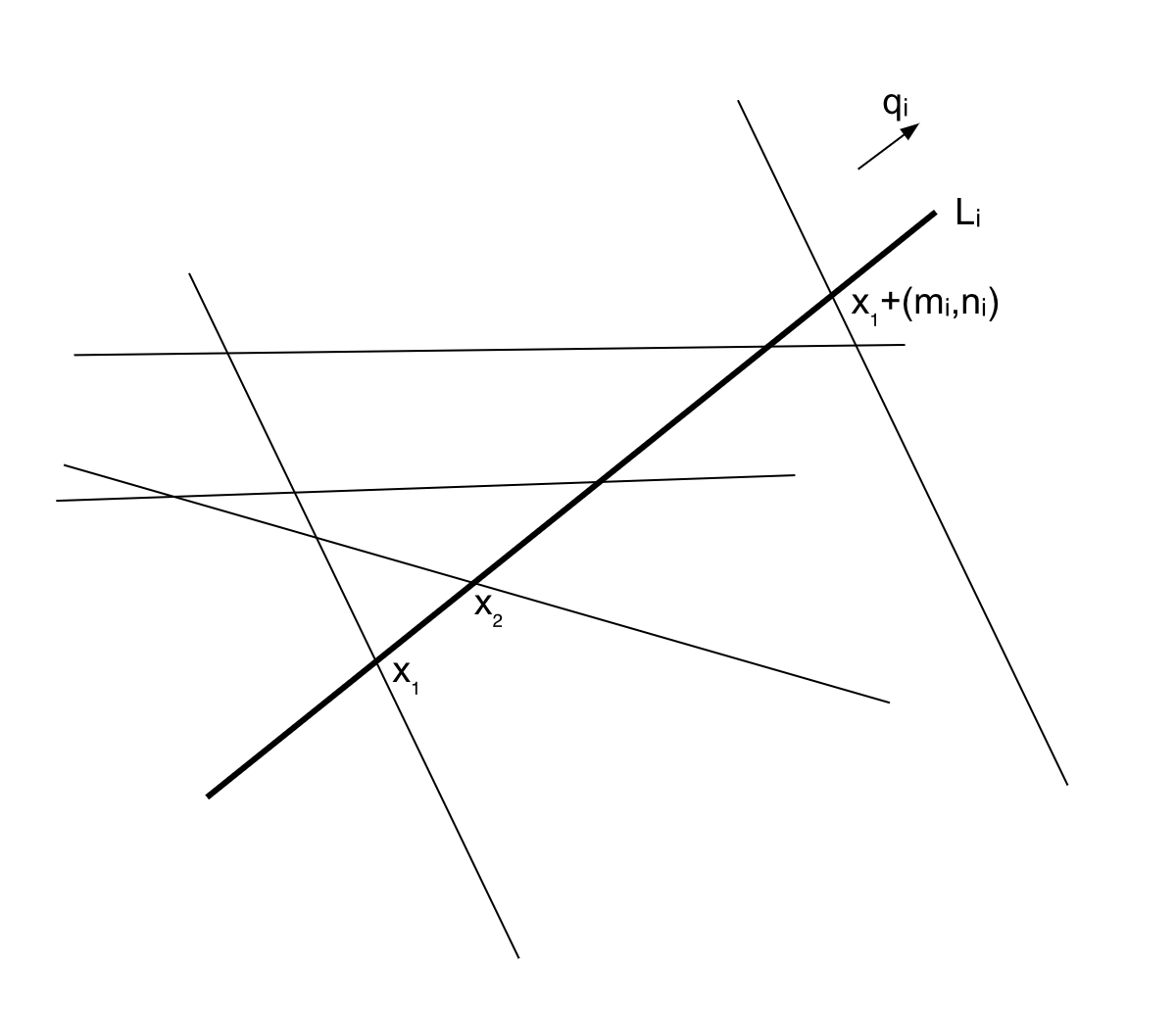}
\captionof{figure}{Intersection points on $L_i$} \label{fig2}
\end{center}

Next, in a small neighborhood of each fixed intersection point in $I$, we perturb the two corresponding intersecting lines a bit to create  gradient flows of an appropriate function. 
Since this is purely local, by linear transformations and translations,  it suffices to show how to perform this procedure in a neighborhood of the origin $(0,0)$ provided that $L_1, L_2$ are the $x_1$-axis and $x_2$-axis, respectively.

Let $\alpha\in (0,1)$ as fixed in Theorem \ref{thm:main}, pick  $k\in \N$ so that
\[
\alpha \leq 1-{1\over 2k}.
\] 
Consider the potential function
\[
u(x_1,x_2)=C_k\left({{x_{1}^{4k}\over C_k}+ x_{2}^{2}}\right)^{1-{1\over 4k}}+2x_1,
\]
where $C_k>1$ is a positive constant to be determined.  
Clearly, $u\in C^{1, 1-{1\over 2k}}(\R^2)$ and  is $C^{2}$ away from the origin.  
\begin{lem} 
Fix $C_k>2k(4k+1)$.
Then, $u$ has infinitely many distinct gradient flows  passing through the origin.
\end{lem}

\begin{proof}
Apparently,  $\gamma_1(t)=(f(t), 0)$ with 
\[
\begin{cases}
f'(t)=2+C_k^{\frac{1}{4k}}(4k-1)f(t)^{4k-2},\\
f(0)=0.
\end{cases}
\]
is a gradient flow of $u$ passing through the origin.

To prove the lemma,  it suffices to show that if $\xi(t)=(x_1(t), x_2(t)):\R\to \R^2$ is a gradient flow of $u$ and 
\[
\xi(0)\in D:=\{(a,b)\,:\, 0<a<1,\  0<b<a^{2k}\}, 
\]
then
\[
\xi((-\infty,0))\cap (0,\infty)^2 \subset D. 
\]
Note that  $x_1(t)$ and $x_2(t)$ are both increasing within $D$ and  $\xi$ cannot intersect with $\gamma_1$ away from the origin. So if the statement is not correct, there exists $\theta<0$ such that 
\[
0<x_2(\theta)=x_{1}^{2k}(\theta) \quad \mathrm{ and } \quad 0<x_2(t)<x_{1}^{2k}(t)<1 \quad \text{for $t\in (\theta, 0)$}.
\]
At $\theta$,
\[
{C_k x_{1}^{2k-1}(\theta)\over 1+4k}<{u_{x_2}(x_1(\theta),x_2(\theta))\over u_{x_1}(x_1(\theta),x_2(\theta))}={x_{2}'(\theta)\over x_{1}'(\theta)}\leq 2 k x_{1}^{2k-1}(\theta).
\]
This contradicts the choice of $C_k$.
The proof is complete.
\begin{center}
\includegraphics[scale=0.7]{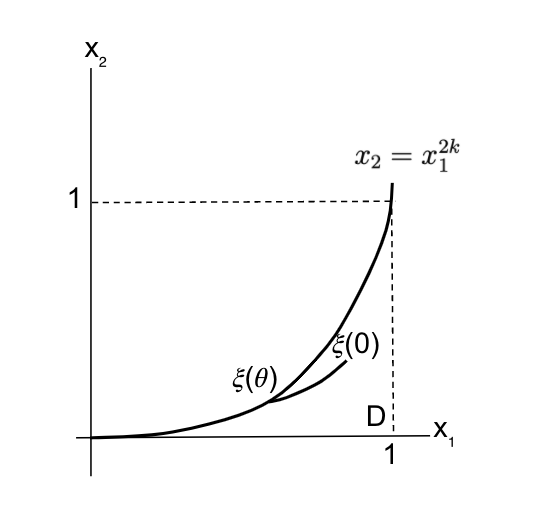}
\captionof{figure}{Graph of $\xi$ in $D$} \label{fig-xi}
\end{center}
\end{proof}


By this construction, we are able to form $m$ periodic  curves $\{\tilde L_i\}_{i=1}^{m}$ and their integer translations such that, for some small  $r \in (0,{d\over 10})$,

\begin{enumerate}

\item $\tilde L_i=L_i$  away from the set $I_r=\{x\in \R^2\,:\,  d(x,I)\leq r\}$;

\item the set of intersection points remains the same, i.e.,  for $i\not= j$ and any integer vector $v\in \Z^2$, 
\[
\tilde L_i\cap (\tilde L_j+v)=L_i\cap (L_j+v);
\]
Equivalently,  $\tilde L_i\cap \tilde L_j=L_i\cap L_j$  when projected to $\Bbb T^2$.

\item given  $i\not= j$ and an integer vector $v\in \Z^2$,  if $\tilde L_i$ and $\tilde L_j+v$ intersect at $x=x_{i,j,v}$, then there exists a $C^{1, \al}$ function $u=u_{i,j,v}$ in $B_{r\over 2}(x)$ such that 

\begin{itemize}
\item $|Du(x)|\geq 1$ in $B_{r\over 2}(x)$;

\item within $B_{r\over 2}(x)$,   $\tilde L_i$ and $\tilde L_j+v$ are two gradient flows of $u$ that only intersect at $x$;

\item(periodicity) if two intersection points $x_{i,j,v}=x_{i',j',v'}+w$ for some $w \in \Z^2$, then
\[
u_{i,j,v}(x+w)=u_{i',j',v'}(x) \quad \text{for $x\in B_r(x_{i',j',v'})$}. 
\]
This says that $u$ is well defined on $I_{r\over 2}$ when being projected to the flat torus $\Bbb T^2$. 
\end{itemize}
\end{enumerate}
The perturbed network is henceforth denoted by
\[
\Gamma=\bigcup_{1 \leq i\leq m}({\tilde L_i+\Z^2}).
\]

\begin{center}
\includegraphics[scale=0.6]{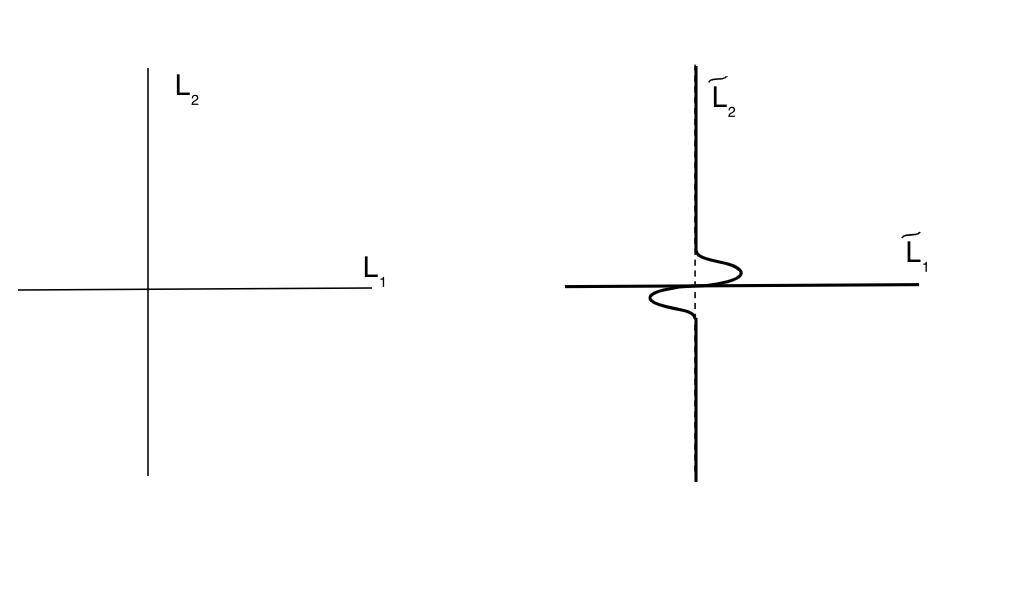}
\captionof{figure}{Local perturbation at the intersection of $L_1$ and $L_2$} \label{fig3}
\end{center}
\smallskip

\noindent {\bf Step 2. Initial choice of $a_0$.}
We can choose $r_0\in (0,{r\over 2})$ and  $a_0 \in C^{1,\alpha}(\Bbb T^2, (0,\infty))$ such that $a_0$ is $C^{\infty}$ away from the set $I$ and satisfies the following conditions. 
\begin{enumerate}

\item for each given intersection point $x=x_{i,j,v}\in I$ and the associated function $u=u_{i,j,v}$ from the above
\[
a_0(y)={1\over |Du(y)|} \quad \text{ for $x\in B_{r_0}(x)$};
\]

\item for every two intersection points $x, y$ on $\tilde L_i$ for $1\leq i\leq m$ (i.e., $x,y\in \tilde L_i\cap I$), the weighted length $l_i(x,y)$ between $x$ and $y$ along $\tilde L_i$ satisfies
\begin{equation}\label{eq:cost}
l_i(x,y) :=\int_{0}^{1}{1\over a_0(\xi(t))}|\dot \xi(t)|\,dt=|p_i\cdot (x-y)|.
\end{equation}
\end{enumerate}
Here, $\xi:[0,1]\to \tilde L_i$ is an arbitrary parametrization of $\tilde L_i$ between $x$ and $y$.  In particular, the weighted length of each period (i.e., from $x$ to $x+(m_i,n_i)$) of  $\tilde L_i$ is $1\over \lambda_i$.  
The existence of $a_0$ is clear provided $r>0$ is small enough. 
By Lemma \ref{one-sided},  
\begin{equation}
\label{eq:lowbdd}
\overline H_{a_0}(p)\geq \max_{1\leq i\leq m}|q_i\cdot p|.
\end{equation}
%
For $i=1,2,...,m$,   let $\xi_i:\R\to \tilde L_i$ be the smooth  reparametrization of $\tilde L_i$ such that 
\[
|\dot \xi_i(t)|={1\over a_0(\xi_i(t))} \quad \text{ for } t\in \R.
\]

Owing to Lemma \ref{edge-construction-1} and the periodicity of $\Gamma$,  there exists a universal  $\delta_0\in (0, r_0)$ such that for each $i=1,2,3,...,m$, there exists $w_i\in C^{1,\alpha}(\tilde L_{i,\delta_0})$ such that $w_i$ is $C^{\infty}$ away from intersection points and 

\begin{enumerate}

\item $\dot \xi_i(t)=Dw_i(\xi_i(t))$ for all $t\in \R$, i.e., $\xi_i$ is the gradient flow  of $w_i$;

\item $Dw_i(x)=Du_{i,j,v}(x)$ for $x\in B_{\delta_0}(x_{i,j,v})$;

\item $\inf_{x\in \tilde L_{i,\delta_0}}|Dw_i(x)|>0$.

\end{enumerate}

\noindent Here, $\tilde L_{i,\delta_0}=\{x \,:\, d(x,\tilde L_i)<\delta_0\}$ and $x_{i,j,v}$ is any intersection point on $\tilde L_i$. 
Let 
\[
\Gamma_{\delta_0}=\{x\in \R^2\,:\, d(x, \Gamma)<\delta_0\} = \bigcup_{i=1}^{m}\,(\tilde L_{i,\delta_0}+\Z^2). 
\]
Then, for $x\in \Gamma_{\delta_0}$,  we define
\[
a_0(x)={1\over |Dw_i(x-v)|} \quad \text{ if $x-v\in \tilde L_{i,\delta_0}$ for $1\leq i \leq m$, and $v \in \Z^2$}.
\]
Extend $a_0$ to $C^{1,\alpha}(\Bbb T^2, (0,\infty))$  in such a way that it is smooth away from $I$. 
\begin{center}
\includegraphics[scale=0.7]{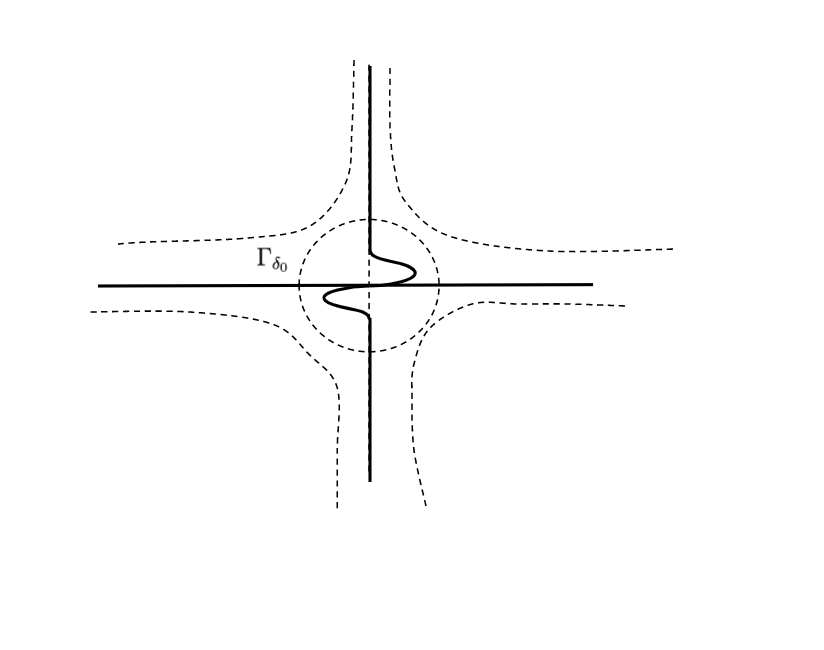}
\captionof{figure}{Part of $\Gam_{\del_0}$} \label{fig-Gam}
\end{center}


\medskip

\noindent {\bf Step 3. Adjustments of $a_0$.} 
Next we need to construct $\tilde a\in C^{1, \alpha}(\Bbb T^2, (0,\infty))$ that is smooth away from $I$,
\[
\tilde a=a_0 \quad \text{ on $\Gamma$},
\]
and,  for $1 \leq i \leq m$,
\[
\overline H_{\tilde a}(p_i)\leq 1.
\]
Since $\tilde a$ agree with $a_0$ of the previous section along $\tilde L_i$'s, the property \eqref{eq:cost} and, by Lemma \ref{one-sided}, the inequality \eqref{eq:lowbdd} are preserved. We hence obtain 
\[
\begin{cases}
\overline H_{\tilde a}(p)\geq \max_{1\leq i\leq m}|q_i\cdot p|,\\[3mm]
\overline H_{\tilde a}(p_i)\leq 1  \quad \text{ for $i=1,2,..., m$}.
\end{cases}
\]
Therefore,  the function $\tilde a$ is exactly what we are looking for, that is, $S_{\tilde a}=P$.

\medskip

Note that, owing to \eqref{strict}, for given $i\in \{1,2,3,..,m\}$, the following points hold

\begin{itemize}
\item for $j=i,i+1$ and two intersection points $x,y\in L_j$,
\[
|p_i\cdot x-p_i\cdot y|=l_j(x,y);
\]

\item  for $j\not=i,i+1$ and  every two distinct intersection points $x,y\in L_j$, 
\[
|p_i\cdot x-p_i\cdot y|=|p_i\cdot (x-y)|\leq \max_{l\not =j-1, j}|p_l\cdot (x-y)|<|p_j\cdot (x-y)|=l_j(x,y).
\]
\end{itemize}
\noindent In light of Remark \ref{edge-construction-2} and the periodicity of $\Gamma$,  there exists $\mu_0\in (0, \delta_0)$ such that  for each $i=1,2,3,..,m$, there exists a function $\tilde u_i\in C^{1, \alpha}(\Gamma_{\mu_0})$ such that 

\[
\begin{cases}
\tilde u_i\in C^{1, \alpha}(\Gamma_{\delta_0}),\ \tilde u_i\in C^{\infty}(\Gamma_{\mu_0}\backslash I), \\
\inf_{\Gamma_{\delta_0}}|D\tilde u_i|>0,\\
\tilde u_i-p_i\cdot x  \quad \text{ is $\Z^2$-periodic in $\Gamma_{\mu_0}$},\\
|D\tilde u_i|\leq |Dw_i|  \quad \text{ in $\Gamma_{\mu_0}$},
\end{cases}
\]
and for any intersection point $x=x_{i,j,v}\in I$,
\[
D\tilde u_i=Dw_i=Du_{i,j,v}  \quad \text{ in $B_{\mu_0}(x_{i,j,v})$}.
\]
We extend $\tilde u_i-p_i\cdot x$ to  $v_i\in C^{1,\alpha}(\Bbb T^2)$ such that $v_i$ is $C^2$ away from $I$, and  for $u_i=p_i\cdot x+v_i$,
\[
u_i=\tilde u_i \quad \text{ on $\Gamma_{\mu_0\over 2}$}.
\]
Now let
\[
K_1=\max_{1\leq i \leq m}\max_{x\in \R^2}|Du_i(x)| \quad \text{ and } \quad K_2=\max_{x\in \R^2}a_0(x).
\]
Choose  $\phi(x)\in C^{\infty}(\Bbb T^2, (0,1])$ such that
\[
\phi(x)=
\begin{cases}
1  \quad &\text{ for $x\in \Gamma_{\mu_0\over 4}$},\\[3mm]
{1\over K_1(1+K_2)} \quad &\text{ for $x\in \R^2\backslash \Gamma_{\mu_0\over 2}$}.
\end{cases}
\]
Finally, let
\[
\tilde a(x)={\phi(x)a_0(x)} \quad \text{ for } x\in \R^2.
\]
Then, for $i=1,2,...,m$,
\[
\begin{cases}
\tilde a(x)|p+Dv_i(x)|\leq \tilde a(x)|Dw_i(x)|=\phi(x)\leq 1 \quad &\text{ for $x\in \Gamma_{\mu_0\over 2}$,}\\[3mm]
\tilde a(x)|p+Dv_i(x)|={a_0(x)|Du_i(x)|\over K_1(1+K_2)}\leq 1 \quad &\text{ for $x\in \R^2\backslash \Gamma_{\mu_0\over 2}$,}
\end{cases}
\]
which says
\[
\max_{x\in \R^2}\,\tilde a(x)|p+Dv_i(x)|=\max_{x\in \R^2}\,\tilde a(x)|Du_i(x)|\leq 1.
\]
\medskip
By the inf-max formula \eqref{inf-max}, for $1 \leq i \leq m$,
\[
\overline H_{\tilde a}(p_i)\leq 1.
\]
This verifies that $\tilde a$ constructed above has the desired properties, and the proof of Theorem \ref{thm:main} is completed.

\begin{rem}\label{newproof}
Our method also provides a simple proof of the following result in \cite{MPS}:
there exists $a\in C^{\infty}(\T^2, (0,\infty))$ such that 
\begin{equation}\label{eq:MPS}
\{\pm q_i\,:\, 1\leq i \leq m\}\subset \partial \ol B_1^a. 
\end{equation}
In fact, to obtain this claim, no gradient matching is needed at the intersections. 
Steps 1--2 in the above proof are not needed. 
Below we give some adaptions to get \eqref{eq:MPS}.
We use the straight line network $\cup_{i=1}^{m}\left(L_i+\Z^2\right)$ directly.

\begin{enumerate}
\item  Pick $a \in C^\infty(\Bbb T^2, (0, \infty))$ such that  $a=1$ in a small neighborhood of $I$, and \eqref{eq:cost} holds with $a, L_i$ in place of $a_0, \tilde L_i$, respectively.

\item In Step 3 of the proof above, choose $u_i$ as  
\[
u_i(x)={q_i\over |q_i|} \cdot (x-x_{i,j,v})+p_i\cdot x_{i,j,v}
\]
 near each intersection point  $x_{i,j,v}$. 
Then, using the method of characteristics (see \cite[Chapter 3]{Ev2} for instance), we extend $u_i(x)$ to a smooth function on $\Gamma_\delta$ for some $\delta>0$ such that 
\[
a(x)|Du_i(x)|=1  \quad \text{ in $\Gamma_\delta$}.
\]
\item Finally, we follow the same arguments as those in Step 3 of the proof above to conclude. 
\end{enumerate}

\end{rem}

\appendix
\section{Some auxiliary lemmas} \label{appen}

\begin{lem}\label{edge-construction-1}
Suppose that $\gamma:[0,1]\to \R^2$ is a smooth curve satisfying that 
\begin{enumerate}
\item $\min_{t\in [0,1]}|\dot \gamma(t)|>0$ and $\gamma(t_1)\not=\gamma (t_2)$ for $t_1\not =t_2$;

\item there exist $r>0$ and $u_0, u_1 \in C^{\infty}(\R^2)$  such that
\[
\begin{cases}
\dot \gamma(t)=Du_0(\gamma(t)) \quad &\text{ for $t\in [0,r]$},\\
\dot \gamma(t)=Du_1(\gamma(t)) \quad &\text{ for $t\in [1-r,0]$}.
\end{cases}
\]
\end{enumerate}
Then, there exist $\delta>0$, an open neighborhood $U$ of $\gamma$, and $u\in C^{\infty}(U)$ such that 
\[
\begin{cases}
\inf_{U}|Du|>0,  \\  
Du=Du_0 \quad &\text{ in $B_{\delta}(\gam(0))$}, \\ 
 Du=Du_1 \quad &\text{ in $B_{\delta}(\gam(1))$}
\end{cases}
\]
and
\[
\dot \gam(t)=Du(\gam (t)) \quad \text{ for $t\in [0,1]$}.
\]
\end{lem}

The proof of the above lemma is standard, and we leave it as an exercise for the interested readers.

\begin{rem}\label{edge-construction-2} 
Consider the same set-up of Lemma \ref{edge-construction-1}.
Let $a(\gam(t))={1\over |Du(\gam(t))|}$ for $t\in [0,1]$, and 
\[
M=\int_{0}^{1}{1\over a(\gam(t))}|\dot \gam(t)|\,dt=u(\gam(1))-u(\gam(0)).
\]
For each $r\in (-M,M)$, let $\tau>0$ be sufficiently small, and choose $h\in C^{\infty}(\R)$ so that 
\[
\begin{cases}
h(t)=t  \quad &\text{ for $t\in [0, {\tau\over 2}]$},\\
h(t)=r+t-M \quad &\text{ for $t\in [M-{\tau\over 2}, M]$},\\
|h'(t)|\leq 1 \quad &\text{ for all $t\in [0, M]$}.
\end{cases}
\]
Then, $u_r=h(u-u(\gam(0)))+u(\gam(0))$ satisfies that 
\[
u_r(\gam(0))=u(\gam(0)),  \quad u_r(\gam(1))=u(\gam(0))+r.
\]
Moreover, we also have $|Du_r| \le |Du|$ in $U$ and
\[
Du_r(x)=Du(x) \quad \text{ for $x\in B_\mu(\gam(0))\cup B_\mu(\gam(1))$}
\]
for some $\mu>0$. 
\end{rem}

\end{proof}


\end{document}